\documentclass{amsart}
\usepackage{amssymb,amsmath,amsfonts,euscript,mathrsfs,amsthm}
\usepackage[dvips]{graphicx}

\usepackage[dvips]{color}

\input xy
\xyoption{all} \CompileMatrices

\DeclareMathOperator{\Id}{Id} \DeclareMathOperator{\lie}{Lie}

\def\cC{{\mathscr C}}

\def\cH{{\mathscr H}}

\def\cT{{\mathscr T}}
\def\cR{{\mathscr R}}

\def\fl{{\mathfrak l}}

\def\fs{{\mathfrak s}}

\def\End{\operatorname{End}}

\def\Spec{\operatorname{Spec}}

\def\res{\operatorname{res}}

\def\vac{|0 \rangle}

\newtheorem{thm}{Theorem}[section]
\numberwithin{equation}{thm}

\newtheorem{lem}[thm]{Lemma}

\theoremstyle{definition}
\newtheorem{defn}[thm]{Definition}

\newtheorem{ex}[thm]{Example}

\theoremstyle{remark}

\newtheorem{rem}[thm]{Remark}

\newtheorem{nolabel}[thm]{}

\def\on{\operatorname}
\begin{document}
\title{Supersymmetry of The chiral de Rham complex}
\author{David Ben-Zvi}
\address{Department of Mathematics, University of Texas,  Austin, TX 78712, USA}
\email{benzvi@math.utexas.edu}
\author{Reimundo Heluani}
\address{Department of Mathematics, MIT, Cambridge, MA 02139, USA}
\email{heluani@math.mit.edu}
\author{Matthew Szczesny}
\address{Department of Mathematics and Statistics, Boston University, Boston, MA 02215, USA}
\email{szczesny@math.bu.edu}

\maketitle
\begin{abstract} We present a superfield formulation of the chiral de Rham complex
(CDR) \cite{malikov} in the setting of a general smooth manifold,
and use it to endow CDR with superconformal structures of geometric
origin. Given a Riemannian metric, we construct an $N=1$ structure
on CDR (action of the $N=1$ super--Virasoro, or Neveu--Schwarz,
algebra). If the metric is K\"{a}hler, and the manifold Ricci-flat,
this is augmented to an $N=2$ structure. Finally, if the manifold is
hyperk\"{a}hler, we obtain an $N=4$ structure. The superconformal
structures are constructed directly from the Levi-Civita connection.
These structures provide an analog for CDR of the extended
supersymmetries of nonlinear $\sigma$--models.
\end{abstract}
\section{Introduction}\label{sec:1}

In the paper \cite{malikov}, Malikov, Schechtman and Vaintrob
introduced a sheaf of vertex superalgebras $\Omega^{ch}_M$ attached
to any smooth complex variety $M$, called the \emph{chiral de Rham
complex} of $M$ (see also \cite{GMSII}). If $M$ is $n$--dimensional,
the fibers of $\Omega^{ch}_M$ are isomorphic as vertex superalgebras
to a completion of the $bc-\beta \gamma$ system on $n$ generators,
or in physics terminology, to the tensor product of the bosonic and
fermionic ghost systems on $n$ generators. The sheaf cohomology of
$\Omega^{ch}_M$, $H^{*}(M,\Omega^{ch}_M)$, also a vertex
superalgebra, is related to the chiral algebra of the half-twist of
the $\sigma$-model with target $M$, a quantum field theory
associated to $M$ (see \cite{FL}, \cite{kapustin}, \cite{Witten} ).
It is shown in \cite{malikov} that in the holomorphic setting, for
arbitrary $M$, $H^{*}(M, \Omega^{ch}_M)$ carries a conformal
structure, and when $M$ has a global holomorphic volume form,
$H^{*}(M,\Omega^{ch}_M)$ admits $N=2$ superconformal symmetry
(equivalently, admits the structure of a topological vertex
algebra).

In this paper, we present a {\em superfield} construction of the
chiral de Rham complex in the $C^{\infty}$ setting, and examine how
various geometric structures on a manifold give rise to extra
symmetries on $\Omega^{ch}_M$. (For the $C^\infty$ formulation of
CDR see \cite{malikov} and \cite{LL} -- we present a streamlined
formulation in Section \ref{Cinfty}. For a related super--spacetime
approach to CDR see \cite{gerbesIII}.) We show that a Riemannian
metric on $M$ gives rise to an $N=1$ structure on $\Omega^{ch}_M$,
i.e. there exist global sections of $\Omega^{ch}_M$ that generate an
$N=1$ superconformal vertex subalgebra. The $N=1$ structure is
constructed from the Levi-Civita connection on $M$. When $M$ is
Ricci-flat and K\"{a}hler, this is augmented to an $N=2$ structure (
quasiclassical limits of the $N=1$ and $N=2$ structures on CDR were
independently obtained by  F. Malikov \cite{private}, starting from
the Lagrangian for the sigma-model). In this case, choosing
holomorphic and anti-holomorphic coordinates on $M$, $\Omega^{ch}_M$
can be locally written as a tensor product of a holomorphic and
anti-holomorphic part. The $N=2$ structure splits as a tensor
product of two commuting $N=2$ vertex algebras, one ``holomorphic''
and another ``anti-holomorphic''. These structures are slightly
different than the ones considered in \cite{malikov}, and agree only
in the case when the metric $g$ is flat.

 Finally, we consider the case where $M$ is a hyperk\"{a}hler
manifold. Recall that a hyperk\"{a}hler manifold is a Riemannian
manifold possessing three isometric complex structures $I, J, K$,
parallel with respect to the Levi-Civita connection, and satisfying
the quaternionic relations:
\[
 I^2 = J^2 = K^2 = - \Id, \; \; IJ=-JI=K
\]
The real dimension of $M$ is then necessarily a multiple of $4$, and
the above is equivalent to the holonomy group of $M$ being contained
in $Sp(n,\mathbb{H}) \subset SO(4n,\mathbb{R})$. We show that in
this case, the sheaf $\Omega^{ch}_M$ carries $N=4$ superconformal
symmetry, i.e., has an embedding of the $N=4$ superconformal vertex
algebra.

 Our approach relies heavily on the superfield formalism
for vertex algebras introduced in \cite{heluani} under the moniker
of SUSY vertex algebras. In the standard approach to vertex
superalgebras one considers \emph{fields} which are
endomorphism-valued distributions
\[
A(z) = \sum_{n\in \mathbb{Z}} z^{-n-1} A_{(n)}  \in \on{End} V [[z, z^{-1}]], \; \;
A_{(n)} \in \on{End} V
\]
where $V$ is the vertex superalgebra. Let $\theta$ be an odd formal variable
satisfying $\theta^2 = 0$. A superfield is an endomorphism-valued distribution of the form
\[
A^{s} (z, \theta) = A(z) + \theta B(z), \; \; A(z), B(z) \in \on{End} V [[z,z^{-1} ]]
\]
The operator products of the two fields $A(z), B(z)$ are now encoded
in the operator product of the single superfield $A^{s}(z,\theta)$.
In the standard field approach to vertex superalgebras, the $N=4$
superconformal algebra is generated by eight fields. Thus,
constructing a representation of this object in terms of free fields
involves checking that a very large number of operator products are
correct. In the SUSY vertex algebra formalism, the $N=4$
superconformal algebra is generated by only four fields, which
greatly simplifies the computations involved.

The SUSY vertex algebra formalism also yields a simpler, tensorial
description of the chiral de Rham complex. In a coordinate patch $U$
diffeomorphic to $\mathbb{R}^n$, $\Omega^{ch}_X (U)$ is generated by
sections $a^i(z), b^i(z), \phi^i(z), \psi^i(z ), \; i=1 \cdots n$.
Under a change of coordinates, $b^i(z)$ transforms as a function,
$\phi^i(z)$ as a 1-form, $\psi^i(z)$ as a vector field, while
$a^i(z)$ transforms in a seemingly complicated, non--linear way (see
\ref{eq:3.1.2}). In the SUSY formalism, these are combined into two
superfields
\[
B^i(z,\theta) = b^i(z) + \theta \phi^i(z), \; \; \; \Psi^i(z,\theta) = \psi^i(z) +
\theta a^i(z)
\]
We show that under a change of coordinates, $B^i(z,\theta)$
transforms as a function and $\Psi^i(z,\theta)$ simply as a vector
field. In particular the well known cancelation of anomalies
necessary for the construction of the chiral de Rham complex becomes
an obvious consequence of the superfield formalism.

To construct the desired supersymmetries, we first show that
endomorphisms of the tangent bundle give sections of the chiral de
Rham complex, and hence (under the state--field correspondence)
fields or superfields. The basic $N=1$ structure (Neveu--Schwarz
current) of CDR on any Riemannian manifold is constructed explicitly
out of the metric. Finally, to enhance this to an $N=2$
(respectively, $N=4$) structure on a Calabi--Yau (respectively,
hyperk\"ahler) manifold we simply adjoin the superfield associated
to the complex structure endomorphism $I$ (respectively, the
superfields associated to $I,J,K$).

The superconformal structures exhibited in this paper are inspired
by the well--known supersymmetries of $\sigma$--models. Zumino
\cite{Zumino} and Alvarez-Gaum\'e and Freedman \cite{AGF} showed
that nonlinear $\sigma$--models on general Riemannian, complex and
hyperk\"ahler targets carry $N=1,2,4$ super--Poincar\'e symmetries,
respectively (see \cite{HKLR,Dan,IAS} for excellent discussions, and
\cite{BGL} for a related recent development). In the $N=2$ case,
this structure is super{\em conformal} when the target is
Calabi-Yau. (Note that the quantum $\sigma$--model is only
conformally invariant, up to higher order terms, when the Ricci
curvature vanishes.)

It is interesting to compare the superconformal structures we
exhibit with the bundles of superconformal vertex algebras
constructed by Tamanoi \cite{Tamanoi} and Zhou \cite{zhou}. Note
that these bundles of vertex algebras are {\em linear}, i.e.
associated to the (frame bundles of the) tangent bundles of
Riemannian manifolds, and should be compared with the associated
graded of our construction with respect to a natural filtration on
CDR (which depends on all order jets on the manifold). Note also
that \cite{malikov} endow CDR with an $N=2$ structure that is
defined for any smooth manifold, coming from the choice of a volume
form on that manifold. However this structure is different from ours
(and that for the $\sigma$--model), which depends directly on the
metric, and does not seem to allow for an $N=4$ extension in the
hyperk\"ahler case.

Finally, the superfield formulation of the chiral de Rham complex is
closely related to its geometric formulation as a {\em factorization
algebra} over super--curves, extending the factorization structure
over (even) curves described by Kapranov and Vasserot \cite{KV}.
This geometric construction was an initial motivation for this work
and we plan to return to it in the future.

The outline of the paper is as follows. In Section \ref{classical}
we recall the standard approach to vertex superalgebras following
\cite{kac:vertex}, and review the $N=1,2,4$ superconformal vertex
algebras in this language. In Section \ref{sec:2} we introduce the
formalism of SUSY vertex algebras following \cite{heluani}, and give
several examples, including a description of the $N=1,2,4$
superconformal vertex algebras in this more compact language.
Section \ref{sec:3} recalls the construction of the chiral de Rham
complex and gives its simplified superfield description. Section
\ref{Cinfty} presents a coordinate--free construction of CDR on
general smooth manifolds. Section \ref{hyperkaler} is a brief review
of hyperk\"{a}hler manifolds. Finally, Section \ref{geometry}
contains the construction of the $N=1,2,4$ structures on the chiral
de Rham complex.

{\bf Acknowledgements}. DB would like to thank the Aspen Center for
Physics for its hospitality. RH would like to thank his advisor
Victor Kac. MS would like to thank Edward Frenkel
and Vassily Gorbounov for valuable discussions, and RIMS for their
hospitality. In June 2005,  MS attended a conference at the Erwin Schrodinger Institute, where Fyodor Malikov explained that he had obtained quasi-classical
limits of the $N=1$ and $N=2$ structures on CDR. MS is very grateful for this, as well as other valuable discussions.
During this project, DB was supported by a Summer
Research Assignment from the University of Texas and MS was
supported by NSF grant DMS-0401619.

\section{Vertex superalgebras} \label{classical}

In this section, we review the definition of vertex superalgebras,
as presented in \cite{kac:vertex}, in order to fix notation, and
facilitate comparison with the super-field formalism introduced
later.

\begin{nolabel}
    Given a vector space $V$, an \emph{$\End(V)$-valued field} is a formal
    distribution of the form
    \begin{equation}
        A(z) = \sum_{n \in \mathbb{Z}} z^{-1-n} A_{(n)},\qquad A_{(n)} \in
        \End(V),
    \end{equation}
    such that for every $v \in V$, we have $A_{(n)}v = 0$ for large enough $n$.
    \label{no:2.1.1a}
\end{nolabel}
\begin{defn}
    A vertex super-algebra consists of the data of a super vector space $V$,
    an even vector $\vac \in V$ (the vacuum vector),
     an even endomorphism $T $, and a parity preserving linear map $A \mapsto Y(A,z)$ from
     $V$ to $\End(V)$-valued fields (the state-field correspondence). This
     data should satisfy the following set of axioms:
    \begin{itemize}
    \item Vacuum axioms:
        \begin{equation}
            \begin{aligned}
            Y(\vac, z) &= \Id \\
            Y(A, z) \vac &= A + O(z) \\
            T \vac &= 0
            \end{aligned}
            \label{eq:1.3.1}
        \end{equation}
    \item Translation invariance:
        \begin{equation}
            \begin{aligned}
                {[}T, Y(A,z)] &= \partial_z Y(A,z)
        \end{aligned}
            \label{eq:1.3.2}
        \end{equation}
    \item Locality:
        \begin{equation}
            (z-w)^n [Y(A,z), Y(B,w)] = 0 \qquad n \gg 0
            \label{eq:1.3.3}
        \end{equation}
     \end{itemize}
(The notation $O(z)$ denotes a power series in $z$ without constant
term.)
    \label{defn:1.3}
\end{defn}
\begin{nolabel}
    Given a vertex super-algebra $V$ and a vector $A \in V$, we expand the fields
    \begin{equation}
        Y(A,z) = \sum_{{j \in \mathbb{Z}}} z^{-1-j}
        A_{(j)}
        \label{eq:1.4.1}
    \end{equation}
    and we call the endomorphisms $A_{(j)}$ the \emph{Fourier modes} of
    $Y(a,z)$. Define now the operations:
    \begin{equation}
        \begin{aligned}
            {[}A_\lambda B] &= \sum_{{j \geq 0}}
            \frac{\lambda^{j}}{j!} A_{(j)}B \\
            A B &= A_{(-1)}B
        \end{aligned}
        \label{eq:1.4.2}
    \end{equation}
    The first operation is called the $\lambda$-bracket and the second is
    called the \emph{normally ordered product}.
     The $\lambda$-bracket contains all of the information about the commutators between the Fourier coefficients of fields in $V$. \label{no:1.4}
\end{nolabel}

\subsection{The $N=1$, $N=2$, and $N=4$ superconformal vertex algebras}

In this section we review the standard description of the $N=1,2,4$
superconformal vertex algebras. In section \ref{sec:2}, the same
algebras will be described in the SUSY vertex algebra formalism.

\begin{ex}{\bf The $N=1$ (Neveu-Schwarz)
    superconformal vertex algebra} \label{N1ex}

    The $N=1$ superconformal vertex algebra (\cite{kac:vertex}) of central charge $c$ is generated by two fields: $L(z)$, an even field of conformal weight $2$, and $G(z)$, an odd primary field of conformal weight $\frac{3}{2}$, with the $\lambda$-brackets
\begin{equation}
{[L}_\lambda L] = (T + 2\lambda) L + \frac{c \lambda^3}{12}
    \label{eq:1}
\end{equation}
\[
{[L}_\lambda G] = (T+\frac{3}{2} \lambda) G
\]
\[
{[G}_\lambda G]  = 2L + \frac{c \lambda^2}{3}
\]
$L(z)$ is called the Virasoro field.
\end{ex}

\begin{ex}{\bf The $N=2$ superconformal vertex algebra} \label{N2ex}

The $N=2$ superconformal vertex algebra of central charge $c$ is
generated by the Virasoro field $L(z)$ with $\lambda$-bracket
(\ref{eq:1}), an even primary field $J(z)$ of conformal weight $1$,
and two odd primary fields $G^{\pm}(z)$ of conformal weight
$\frac{3}{2}$, with the $\lambda$-brackets \cite{kac:vertex}
\begin{equation}
    {[}L_\lambda J] = (T + \lambda) J
    \label{eq:2}
\end{equation}
\begin{equation}
    [L_\lambda G^\pm] = \left(  T + \frac{3}{2}\lambda \right) G^\pm
    \label{eq:3}
\end{equation}
\begin{xalignat}{2}
    {[J}_\lambda G^\pm] &= \pm G^\pm & [J_\lambda J] &= \frac{c}{3}\lambda \\
    {[G^+}_\lambda G^-] &= L + \frac{1}{2} TJ + \lambda J + \frac{c}{6}\lambda^2
    & {[G^\pm}_\lambda G^\pm] &= 0
    \label{eq:4}
\end{xalignat}

\end{ex}

\begin{ex}{\bf The ``small'' $N=4$ superconformal vertex algebra} \label{N4ex}

    The even part of this vertex algebra is generated by the Virasoro field
    $L(z)$ and three primary fields of conformal weights $1$, $J^0$, $J^+$ and $J^-$. The odd
    part is generated by four primary fields of conformal weight $3/2$, $G^\pm$
    and $\bar{G}^\pm$. The remaining
(non-vanishing)
$\lambda$-brackets are (cf \cite[page 36]{kacwakimoto1})
\begin{xalignat}{2}
    {[J^0}_\lambda J^\pm] &= \pm 2 J^\pm & {[J^0}_\lambda J^0] &=
    \frac{c}{3}\lambda \\
    {[J^+}_\lambda J^-] &= J^0 + \frac{c}{6} \lambda & {[J^0}_\lambda G^\pm] &= \pm
    G^\pm \\
    {[J^0}_\lambda \bar{G}^\pm] &= \pm \bar{G}^\pm & {[J^+}_\lambda G^-] &= G^+
    \\ {[J^-}_\lambda G^+]&= G^- & {[J^+}_\lambda \bar{G}^-] &= - \bar{G}^+ \\
    {[J^-}_\lambda \bar{G}^+] &= - \bar{G}^- & {[G^\pm}_\lambda \bar{G}^\pm]&=
    (T + 2 \lambda) J^\pm \\
    {[G^\pm}_\lambda \bar{G}^{\mp}]&= L \pm \frac{1}{2}T J^0 \pm \lambda J^0 +
    \frac{c}{6} \lambda^2
    \label{eq:6}
\end{xalignat}
\end{ex}
(Note that the $J$ currents form an $\fs\fl_2$ current
algebra.)

\section{SUSY vertex algebras}\label{sec:2} In this section we collect some
results on SUSY vertex algebras (SUSY VAs) from \cite{heluani}. Since we only need
the case with one odd variable, we will adapt the notation to this case, and avoid the prefix
\emph{super} when possible.
\subsection{Structure theory of SUSY VAs}
\begin{nolabel}
    Let us fix notation first. We introduce formal variables $Z=(z,\theta)$ and $W =
    (w,\zeta)$, where $\theta, \zeta$ are odd
    anti-commuting variables and $z, w$ are even commuting variables.
    Given an integer $j$ and $J = 0$ or $1$ we put $Z^{j|J} = z^j \theta^J$.

    Let $\cH$ be the universal enveloping algebra of the
    $1|1$ dimensional Lie super-algebra $[\chi, \chi] = - 2 \lambda$, where $\chi$ is
    odd and $\lambda$ is even and central (super--Heisenberg or Clifford algebra). We will consider another set of
    generators $-S, -T$ for $\cH$ where $S$ is odd, $T$ is central, and $[S, S]
    = 2 T$. Whenever we treat $\cH$ as an $\cH$-module it will be by the
    adjoint action. Denote $\Lambda = (\lambda, \chi)$,
    $\nabla = (T, S)$, $\Lambda^{j|J} = \lambda^j \chi^J$ and $\nabla^{j|J} =
    T^j S^J$.

    Given a super vector space $V$ and a vector $a \in V$, we will denote by
    $(-1)^a$ its parity.
    \label{no:2.1}
\end{nolabel}
\begin{nolabel}
    Let $U$ be a vector space, a $U$-valued formal distribution is an
    expression of the form
    \begin{equation}
        \sum_{\stackrel{j \in \mathbb{Z}}{J = 0,1}} Z^{-1-j|1-J} w_{(j|J)}
        \qquad w_{(j|J)} \in U.
        \label{eq:2.2.1}
    \end{equation}
    The space of such distributions will be denoted by $U[ [Z, Z^{-1}] ]$. If
    $U$ is a Lie algebra we will say that two such distributions $a(Z), \,
    b(W)$ are
    \emph{local} if
    \begin{equation}
        (z - w)^n [a(Z), b(W)] = 0 \qquad n \gg 0.
        \label{eq:2.2.2}
    \end{equation}
    The space of distributions such that only finitely many negative powers
    of $z$ appear (i.e. $w_{(j|J)} = 0$ for large enough $j$) will be denoted
    $U( (Z ))$. In the case when $U = \End(V)$ for another vector space $V$,
    we will say that a distribution $a(Z)$ is a \emph{field} if $a(Z)v \in V(
    (Z ))$ for all $v \in V$.
    \label{no:2.2}
\end{nolabel}
\begin{defn}
    An $N=1$ SUSY vertex algebra consists of the data of a vector space $V$,
    an even vector $\vac \in V$ (the vacuum vector), an odd endomorphism
    $S$ (whose square is an even endomorphism we denote $T$),
    and a parity preserving linear map $A \mapsto Y(A,Z)$ from
     $V$ to $\End(V)$-valued fields (the state-field correspondence). This
     data should satisfy the following set of axioms:
    \begin{itemize}
    \item Vacuum axioms:
        \begin{equation}
            \begin{aligned}
            Y(\vac, Z) &= \Id \\
            Y(A, Z) \vac &= A + O(Z) \\
            S \vac &= 0
            \end{aligned}
            \label{eq:2.3.1}
        \end{equation}
    \item Translation invariance:
        \begin{equation}
            \begin{aligned}
            {[} S, Y(A,Z)] &= (\partial_\theta - \theta \partial_z)
            Y(A,Z)\\
            {[}T, Y(A,Z)] &= \partial_z Y(A,Z)
        \end{aligned}
            \label{eq:2.3.2}
        \end{equation}
    \item Locality:
        \begin{equation}
            (z-w)^n [Y(A,Z), Y(B,W)] = 0 \qquad n \gg 0
            \label{eq:2.3.3}
        \end{equation}
     \end{itemize}
    \label{defn:2.3}
\end{defn}
\begin{rem}
    Given the vacuum axiom for a SUSY vertex algebra, we will use the state
    field correspondence to identify a vector $A \in V$ with its corresponding
    field $Y(A,Z)$.
    \label{rem:nosenosenose}
\end{rem}
\begin{nolabel}
    Given a $N=1$ SUSY vertex algebra $V$ and a vector $A \in V$, we expand the fields
    \begin{equation}
        Y(A,Z) = \sum_{\stackrel{j \in \mathbb{Z}}{J = 0,1}} Z^{-1-j|1-J}
        A_{(j|J)}
        \label{eq:2.4.1}
    \end{equation}
    and we call the endomorphisms $A_{(j|J)}$ the \emph{Fourier modes} of
    $Y(A,Z)$. Define now the operations:
    \begin{equation}
        \begin{aligned}
            {[}A_\Lambda B] &= \sum_{\stackrel{j \geq 0}{J = 0,1}}
            \frac{\Lambda^{j|J}}{j!} A_{(j|J)}B \\
            A B &= A_{(-1|1)}B
        \end{aligned}
        \label{eq:2.4.2}
    \end{equation}
    The first operation is called the $\Lambda$-bracket and the second is
    called the \emph{normally ordered product}.
    \label{no:2.4}
\end{nolabel}
\begin{rem}
    As in the standard setting, given a SUSY VA $V$ and a vector $A \in V$, we
    have:
    \begin{equation}
        Y(TA, Z) = \partial_z Y(A,Z) = [T, Y(A,Z)]
    \end{equation}
    On the other hand, the action of the ``odd'' derivation $S$ is described
    by:
    \begin{equation}
        Y(SA,Z) = \left( \partial_\theta + \theta \partial_z \right) Y(A,Z)
        \neq [S, Y(A,Z)].
    \end{equation}
    \label{rem:caca5}
\end{rem}
The relation with the standard field formalism is as follows.
Suppose that $V$ is a vertex superalgebra  as defined in section
\ref{classical}, together with a homomorphism from the $N=1$
superconformal vertex algebra in example \ref{N1ex}. $V$ therefore
possesses an even vector $\nu$ of conformal weight $2$, and an odd
vector $\tau$ of conformal weight $\frac{3}{2}$, whose associated
fields
\begin{equation}
    \begin{aligned}
 Y(\nu,z) &= L(z) = \sum_{n \in \mathbb{Z}} L_n z^{-n-2} \\
Y(\tau,z) &= G(z) = \sum_{n \in 1/2 + \mathbb{Z}} G_n z^{-n - \frac{3}{2}}
\end{aligned}
\end{equation}
have the $\lambda$-brackets as in example \ref{N1ex}, and where we
require $G_{-1/2}=S$ and $L_{-1}=T$. We can then endow $V$ with the
structure of an $N=1$ SUSY vertex algebra via the state-field
correspondence \cite{kac:vertex}
\[
Y(A,Z) = Y^{c} (A,z) + \theta Y^{c}(G_{-1/2} A, z)
\]
where we have written $Y^{c}$ to emphasize that this is the
``classical" state-field (rather than state--superfield)
correspondence in the sense of section \ref{classical}.

(Note however that there exist $SUSY$ vertex algebras without such a map
from the $N=1$ superconformal vertex algebra.)
\begin{defn}
    Let $\cH$ be as before. \emph{An $N=1$ SUSY Lie
    conformal  algebra} is a $\cH$-module $\cR$ with an operation
    $[\,_{\Lambda}\,]: \cR \otimes \cR \rightarrow \cH
    \otimes \cR$ of degree
    $1$ satisfying:
    \begin{enumerate}
        \item Sesquilinearity
            \begin{equation}
                [S a_\Lambda b] =  \chi [a_\Lambda b]
                \qquad [a_\Lambda S b] = -(-1)^{a} \left(S
                + \chi
                \right) [a_\Lambda b]
                \label{eq:k.sesqui.1}
            \end{equation}
        \item Skew-Symmetry:
            \begin{equation}
                [b_\Lambda a] =  (-1)^{a b} [b_{-\Lambda -
                \nabla} a]
                \label{eq:k.skew.1}
            \end{equation}
            Here the bracket on the right hand side is computed as
            follows: first compute $[b_{\Gamma}a]$, where $\Gamma =
            (\gamma, \eta)$ are generators of $\cH$ super commuting
            with $\Lambda$, then replace $\Gamma$ by $(-\lambda - T,
            -\chi - S)$.
        \item Jacobi identity:
            \begin{equation}
                [a_\Lambda [b_\Gamma c]] = -(-1)^{a} \left[
                [ a_\Lambda b]_{\Gamma + \Lambda} c \right] +
                (-1)^{(a+1)(b+1)} [b_\Gamma [a_\Lambda c]]
                \label{eq:k.jacobi.1}
            \end{equation}
            where the first bracket on the right hand side is computed as in Skew-Symmetry
            and the identity is an identity in $\cH^{\otimes 2} \otimes\cR$.
    \end{enumerate}
    \label{defn:k.conformal.1}

    Morphisms of SUSY Lie conformal algebras are $\cH$-module morphisms
    $\varphi:\cR \rightarrow \cR'$ such that the following diagram is
    commutative:
    \begin{equation}
\xymatrix{
    \cR \otimes \cR \ar[r]^{\varphi \otimes \varphi} \ar[d] & \cR' \otimes \cR'
    \ar[d] \\
    \cH \otimes \cR \ar[r]_{1 \otimes \varphi} & \cH \otimes \cR'
    }
    \end{equation}
\end{defn}
\begin{rem}
    In this definition we consider $\cR \otimes \cR$ as a module over $\cH$
    using the co-multiplication of $\cH$.  Similarly $\cH
    \otimes \cR$ is a module over $\cH$ (recall that $\cH$ is a module over
    itself with the adjoint action). The bracket ${[}\,_{\Lambda}\,]$ is a
    morphism of $\cH$-modules. The Jacobi identity is an identity in $\cH
    \otimes \cH \otimes \cR$
    \label{rem:2.6}
\end{rem}
\begin{nolabel}
    Given an $N=1$ SUSY VA, it is canonically an $N=1$ SUSY Lie conformal algebra with
    the bracket defined in (\ref{eq:2.4.2}). Moreover, given an $N=1$ Lie
    conformal algebra $\cR$, there exists a unique $N=1$ SUSY VA called the
    \emph{universal enveloping SUSY vertex algebra of $\cR$} with the property
    that if $W$ is another $N=1$ SUSY VA and $\varphi : \cR \rightarrow W$ is a
    morphism of Lie conformal algebras, then $\varphi$ extends uniquely to a
    morphism $\varphi: V \rightarrow W$ of SUSY VAs.
    \label{no:2.7}
\end{nolabel}
\begin{nolabel}
    The operations (\ref{eq:2.4.2}) satisfy:
    \begin{itemize}
        \item Quasi-Commutativity:
            \begin{equation}
                ab - (-1)^{ab} ba = \int_{-\nabla}^0 [a_\Lambda
                b] d\Lambda
                \label{eq:2.8.1}
            \end{equation}
        \item Quasi-Associativity
            \begin{equation}
                (ab)c - a(bc) = \sum_{j \geq 0}
                a_{(-j-2|1)}b_{(j|1)}c + (-1)^{ab} \sum_{j \geq 0}
                b_{(-j-2|1)} a_{(j|1)}c
                \label{eq:2.8.2}
            \end{equation}
        \item Quasi-Leibniz (non-commutative Wick formula)
            \begin{equation}
                [a_\Lambda bc ] = [a_\Lambda b] c + (-1)^{(a+1)b}b
                [a_\Lambda c] + \int_0^\Lambda [ [a_\Lambda
                b]_\Gamma c] d \Gamma
                \label{eq:2.8.3}
            \end{equation}
    \end{itemize}
    where the integral $\int d\Lambda$ is $\partial_\chi \int d\lambda$. In
    addition, the vacuum vector is a unit for the normally ordered product
    and the endomorphisms $S, T$ are odd and even derivations respectively of
    both operations.
    \label{no:2.8}
\end{nolabel}
\subsection{Examples}
\begin{ex}
    Let $\cR$ be the free $\cH$-module generated by an odd vector $H$.
    Consider the following Lie conformal algebra structure in $\cR$:
    \begin{equation}
        {[}H_\Lambda H] = (2T + \chi S + 3 \lambda) H
        \label{eq:2.9.1}
    \end{equation}
    This is the \emph{Neveu-Schwarz} algebra (of central charge 0). This
    algebra admits a central extension of the form:
    \begin{equation}
        {[}H_\Lambda H] = (2T + \chi S + 3\lambda) H + \frac{c}{3} \chi \lambda^2
        \label{eq:2.9.2}
    \end{equation}
    where $c$ is any complex number. The associated universal enveloping
    SUSY VA is the \emph{Neveu-Schwarz} algebra of central charge
    $c$\footnote{Properly speaking, we consider the universal enveloping
    SUSY vertex algebra of $\cR \oplus \mathbb{C} C$ with $C$ central and $TC
    = S C = 0$ and then we quotient by the ideal generated by $C = c$ for
    any complex number $c$}. If we
    decompose the corresponding field
    \begin{equation}
        H(z,\theta) = G(z) + 2 \theta L(z)
        \label{eq:2.9.3}
    \end{equation}
    then the fields $G(z)$ and $L(z)$ satisfy the commutation relations of
    the well known $N=1$ super vertex algebra in example \ref{N1ex}.
    \label{ex:2.9}
\end{ex}
\begin{ex}
    Consider now the free $\cH$ module generated by even vectors $\{B^i\}_{i
    = 1}^n$ and odd vectors $\{\Psi^i\}_{i=1}^n$ where the only non-trivial
    commutation relations are:
    \begin{equation}
        {[}B^i_\Lambda \Psi^j] = \delta_{ij} = [\Psi^j_\Lambda B^i]
        \label{eq:2.10.1}
    \end{equation}
    Expand the corresponding fields as:
    \begin{equation}
        B^i(z,\theta) = b^i(z) + \theta \phi^i(z) \qquad \Psi^i(z,\theta)
        = \psi^i(z) + \theta a^i(z)
        \label{eq:2.10.2}
    \end{equation}
    then the fields $b^i$, $a^i$, $\phi^i$ and $\psi^i$ generate the
    $bc-\beta\gamma$ system as in \cite{malikov}.
    \label{ex:2.10}
\end{ex}
\begin{ex} \label{ex:2.11.a}
    The $N=2$ superconformal vertex algebra is generated by $4$
    fields \cite{kac:vertex}. In this context it is generated by two superfields -- an $N=1$ vector $H$ as in
    \ref{ex:2.9} and an even current $J$, primary of conformal weight $1$,
    that is:
    \begin{equation}
        {[}H_\Lambda J] = (2T + 2\lambda + \chi S) J.
        \label{eq:2.11.a.1}
    \end{equation}
    The remaining commutation relation is
    \begin{equation}
        [J_\Lambda J] = - (H + \frac{c}{3} \lambda \chi).
        \label{eq:2.11.a.2}
    \end{equation}
    Note that given
    the \emph{current} $J$ we can recover the $N=1$ vector $H$.
     In terms of the fields of \ref{N2ex}, $H, J$ decompose as

\begin{equation}
    \begin{aligned}
         J(z,\theta) &= - \sqrt{-1}J(z) - \sqrt{-1}\theta \left( G^-(z) - G^+(z)
        \right) \\
        H(z,\theta) &= \left( G^+(z) + G^-(z) \right) + 2 \theta L(z)
    \end{aligned}
    \label{eq:5}
\end{equation}

\end{ex}
\begin{ex} \label{ex:2.11}
    The ``small'' $N=4$ superconformal vertex algebra is a vertex algebra generated by 8
    fields \cite{kac:vertex}. In this formalism, it is generated by four
    superfields $H, J^i$, $i = 0,1,2$, such that each pair $(H, J^i)$ forms an
    $N=2$ SUSY VA as in the previous example and the remaining commutation
    relations are:
    \begin{equation}
            {[}J^i_\Lambda J^j] =
            \varepsilon^{ijk} (S + 2\chi)J^k \qquad i \neq j
        \label{eq:2.11.1}
    \end{equation}
    where $\varepsilon$ is the totally antisymmetric
    tensor. (In other words, we're writing the $N=4$ algebra in terms of an ${\mathfrak su}_2$ basis $J^i$ of
    superfields
    rather than the $\fs\fl_2$ basis $J^0,J^\pm$ of even fields, together with odd $G$ fields, before.)
     In terms of the fields of example \ref{N4ex}, $H, J^{i}$ decompose as

\begin{equation}
    \begin{aligned}
        J^0(z,\theta) &= - \sqrt{-1} J^0(z) - \sqrt{-1} \theta \left(
        \bar{G}^-(z) - G^+(z)
        \right) \\
        J^1(z,\theta) &= \sqrt{-1} \left( J^+(z) + J^-(z) \right) +
        \sqrt{-1} \left( \bar{G}^+(z) - G^-(z) \right) \\
        J^2(z,\theta) &= \left(J^+(z) - J^-(z)\right) + \theta \left(
        \bar{G}^+(z) + G^-(z)
        \right) \\
        H(z,\theta) &= \left(G^+(z) + \bar{G}^-(z) \right) + 2 \theta L(z)
    \end{aligned}
    \label{eq:7}
\end{equation}

\end{ex}

\section{Chiral de Rham Complex}\label{sec:3}
In this section we recollect some results from \cite{malikov}. We
then provide and discuss a super-field formulation of the chiral de
Rham complex in the algebraic setting. For the applications we have
in mind, we will need to work in the $C^\infty$ setting, which is
described in Section \ref{Cinfty}.
\begin{nolabel}
    The Chiral de Rham complex $\Omega^\mathrm{ch}_M$ is a sheaf of vertex
    algebras defined over any
    smooth algebraic variety $M$ over the complex numbers. In order to
    construct such a sheaf, the authors in \cite{malikov} first construct
    a sheaf of super vertex algebras on $\mathbb{C}^n$ and then show that
    we can \emph{glue} these sheaves by studying the action of changes of
    coordinates.

    To construct the sheaf on $U=\Spec \mathbb{C}[x^1, \dots ,x^n]$, we first look at its global
    sections. This vertex algebra can be described in terms of generators
    and relations as follows.
    $\Omega^\mathrm{ch} (U)$ is then the $bc-\beta\gamma$-system
    vertex algebra. Namely, it is generated by fields $\{a^i, b^i, \psi^i,
    \phi^i\}_{i = 1}^n$, with commutation relations
    \begin{equation}
        [a^i_\lambda b^j] = [\phi^i_\lambda \psi^j] = \delta_{ij}
        \label{eq:3.1.1}
    \end{equation}
    where we have identified the coordinate functions $x^i$ with the
    $(-1)$-Fourier mode of the fields $b^i(z)$ (recall that we identify
    vectors in our vertex algebras with the corresponding fields by the
    state-field correspondence).

    The next step in \cite{malikov} is to consider a localization of this
    vertex algebra, whereby we allow expressions of
    the form
    \begin{equation} \label{expr}
    f(b^1(z), \cdots, b^n(z))
    \end{equation}
     where $f(x^1, \cdots, x^n)$ is an arbitrary algebraic function on
     $U$. This allows us to construct a Zariski sheaf on $\mathbb{C}^n$.
     We may also pass to a formal completion, allowing  $f$ to be an arbitrary
     function on the formal disk $\mathrm{Spf} \mathbb{C}[ [x^1, \dots,
     x^n]]$

     Finally, in order to \emph{glue} these sheaves, one has to analyze
     how these generators transform under changes of coordinates of the
     formal disk.
    Given such a change of coordinates $\tilde{x^i} = g^i(x)$, with inverse $x^i
    = f^i(\tilde{x})$ the generating fields transform as:
    \begin{equation}
        \begin{aligned}
            \tilde{b^i} &= g^i(b)  \\
            \tilde{\phi^i} &= \left( \frac{\partial g^i(b)}{\partial
            b^j} \phi^j\right) \\
            \tilde{\psi^i} &= \left( \frac{\partial f^j}{\partial
            \tilde{b^i}}(g(b)) \psi^j \right) \\
            \tilde{a^i} &= \left( a^j \frac{\partial f^j}{\partial
            \tilde{b^i}}(g(b)) \right) + \left( \frac{\partial^2
            f^k}{\partial \tilde{b^j} \partial \tilde{b^l}} (g(b))
            \frac{\partial g^l}{\partial b^r} \phi^r \psi^k \right)
        \end{aligned}
        \label{eq:3.1.2}
    \end{equation}
    \label{no:3.1}
\end{nolabel}
\begin{nolabel}
    Let us analyze this sheaf of vertex algebras as a sheaf of $N=1$ SUSY VAs.
    For this, we combine the generators into super-fields as:
    \begin{equation}
        B^i = b^i + \theta \phi^i \qquad \Psi^i = \psi^i + \theta a^i
        \label{eq:3.2.1}
    \end{equation}
    These fields generate a SUSY VA as in example \ref{ex:2.10}. Given a
    change of coordinates $g$ as above, the formulas (\ref{eq:3.1.2}) imply that these
    fields transform as:
    \begin{equation}
        \begin{aligned}
        \tilde{B^i} &= g^i(B) \\
        \tilde{\Psi^i} &= \left( \frac{\partial f^j}{\partial
        \tilde{B^i}}(g(B)) \Psi^j \right)
        \end{aligned}
        \label{eq:3.2.4}
    \end{equation}
    Therefore the chiral de Rham complex is described in a simple fashion
    when viewed as a sheaf of SUSY VAs.

    Conversely, the transformation properties (\ref{eq:3.2.4}) imply
    (\ref{eq:3.1.2}). Indeed, suppose the fields $\Psi^i$ and $B^i$ transform
    as in (\ref{eq:3.2.4}). Evaluating at $\theta = 0$ we immediately obtain
    the transformation properties of $\psi^i$ and $b^i$ as in (\ref{eq:3.1.2}).
    We now note that
    \begin{equation}
        S\tilde{B}^i = S g^i(B) = \frac{\partial g^i(B)}{\partial B^j} SB^j
    \end{equation}
    and evaluating at $\theta = 0$ we obtain the transformation property of
    $\phi^i$. Finally, we have:
    \begin{equation}
        \begin{aligned}
        S\tilde\Psi^i &= S \left( \frac{\partial f^j}{\partial
        \tilde{B}^i}(g(B)) \Psi^j \right) \\
        &= \left(\frac{\partial^2 f^j}{\partial \tilde{B}^i \partial\tilde{B}^k}
        (g(B)) S \tilde{B}^k \right) \Psi^j + \frac{\partial f^j}{\partial
        \tilde{B}^i}(g(B)) S\Psi^j\\
        &= \left(\frac{\partial^2 f^j}{\partial \tilde{B}^i \partial\tilde{B}^k}
        (g(B))\frac{ \partial g^k(B)}{\partial B^l} S B^l \right) \Psi^j + \frac{\partial f^j}{\partial
        \tilde{B}^i}(g(B)) S\Psi^j
    \end{aligned}
    \label{eq:caca2}
    \end{equation}
    Using quasi-commutativity we can write the second term as
    \begin{equation}
        S\Psi^j \frac{\partial f^j}{\partial \tilde{B}^i} (g(B)) -
        T\left(\frac{\partial^2 f^j}{\partial \tilde{B}^i \partial \tilde{B}^j}
        \left( g(B) \right) \frac{\partial \tilde{B}^l}{ \partial B^j}
        \right)
    \end{equation}
    On the other hand, using quasi-associativity we can write the first term of
    (\ref{eq:caca2}) as:
    \begin{equation}
        \frac{\partial^2 f^j}{\partial \tilde{B}^i \partial\tilde{B}^k}
        (g(B))\frac{ \partial g^k(B)}{\partial B^l} S B^l  \Psi^j + T \left(\frac{\partial^2 f^j}{\partial \tilde{B}^i \partial\tilde{B}^k}
        (g(B))\frac{ \partial g^k(B)}{\partial B^j} \right)
    \end{equation}
    Adding these two terms and evaluating at $\theta = 0$ we obtain the
    transformation property of $a^i$.

\subsubsection{}
    From this perspective, we can construct the chiral de Rham complex as
    a sheaf of SUSY vertex algebras, by arguing as in \cite{malikov},
    replacing vertex algebras by SUSY VAs and using (\ref{eq:3.2.4})
    instead of (\ref{eq:3.1.2}).  To do this we must check that (\ref{eq:3.2.4}) preserves the SUSY VA
    structure. This can be done as follows:
     we check immediately that
    \begin{equation}
        {[\tilde{B^i}}_\Lambda \tilde{B^j}] = 0
        \label{eq:3.4a.1}
    \end{equation}
    On the other hand, we have from Wick formula:
    \begin{equation}
        {[ \tilde{B^i}}_\Lambda \tilde{\Psi^j}] = \frac{\partial
        f^k}{\partial \tilde{B^j}} ( g(B)) \frac{\partial
        g^i(B)}{\partial(B^k)} = \delta_{ij}
        \label{eq:3.4b.2}
    \end{equation}
    To compute $[\tilde{\Psi^i}_\Lambda \tilde{\Psi^j}]$ we first need
    \begin{equation}
        {[\tilde{\Psi^i}}_\Lambda \frac{\partial f^k}{\tilde{B^j}}
        (g(B))] = \frac{\partial f^l}{\partial \tilde{B^i}} (g(B))
        \frac{\partial}{\partial B^l} \left(  \frac{\partial
        f^k}{\partial \tilde{B^j}} (g(B)) \right)
        \label{eq:3.4b.3}
    \end{equation}
    and
    \begin{equation}
        {[\tilde{\Psi^i}}_\Lambda \Psi^k] = - \frac{\partial}{\partial
        B^k} \left( \frac{\partial f^l}{\partial \tilde{B^i}}(g(B))
        \right) \Psi^l
        \label{eq:3.4b.4}
    \end{equation}
    Now using the Wick formula and noting that the integral term
    vanishes we obtain:
    \begin{equation}
        \begin{aligned}
            {[\tilde{\Psi^i}}_\Lambda \tilde{\Psi^j}] &=
            \frac{\partial f^l}{\partial \tilde{B^i}} (g(B))
            \frac{\partial}{\partial B^l} \left(
            \frac{\partial f^k}{\partial \tilde{B^j}}(g(B)) \right)
            \Psi^k -
            i \leftrightarrow j \\
            &= \frac{\partial f^l}{\partial \tilde{B^i}} (g(B))
            \frac{\partial^2 f^k}{\partial \tilde{B^j} \partial
            \tilde{B^m}} (g(B)) \frac{\partial g^m(B)}{ \partial B^l}
            \Psi^k
            - i \leftrightarrow j \\
            &= \frac{\partial^2 f^k}{\partial \tilde{B^j} \partial
            \tilde{B^i}} (g(B)) \Psi^k - i \leftrightarrow j \\
            &= 0
        \end{aligned}
        \label{eq:3.4b.5}
    \end{equation}
    The equivalence of (\ref{eq:3.1.2}) and (\ref{eq:3.2.4}) shows that the
    sheaf constructed is in fact the chiral de Rham complex of $M$.

\subsubsection{Remark.}     We note that in this
    approach all the \emph{cancellation of anomalies} is contained in the fact
    that the integral term in the $\Lambda$-bracket (\ref{eq:3.4b.5}) vanishes, which in turn is an
    obvious consequence of our formalism.
    \label{rem:3.4b}
    \label{no:3.2}
\end{nolabel}

\subsubsection{Formal setting.}
\begin{nolabel} In the formal setting, the chiral de Rham complex is
    constructed by using the standard arguments of ``formal
    geometry'' \cite[3.9]{malikov}.
    (see also \cite[ch. 17]{frenkelzvi}), i.e. using an action of
    the Lie algebra of vector fields on the formal $n$-dimensional
    disk on the $bc-\beta\gamma$-system. Indeed the vector field $f(x_i)
    \partial_{x_j}$ acts as the residue of the field
    \begin{equation}
        f(b^i) a^j + \sum_{k = 1}^n \left( \partial_{x_k} f \right)
    (b^i)\phi^k \psi^j.
    \end{equation}

    In the context of SUSY vertex algebras, the vector field $f(x_i)
    \partial_{x_j}$ simply acts as the \emph{super residue} of the super-field
    \begin{equation}
        f(B^i) \Psi_j
    \end{equation}
 Here the super residue is defined to be
    \begin{equation}
        \mathrm{sres}_{z,\theta} f(z,\theta) = \partial_\theta \res_z
    f(z,\theta).
    \end{equation}
    \label{no:caca3a}
\end{nolabel}
\section{The $C^\infty$ case} \label{Cinfty}
\begin{nolabel}
    In this section we give a \emph{coordinate independent}
    description of the \emph{chiral de Rham complex} of a smooth
    differentiable manifold $M$. This construction is essentially a
    super-field reformulation of the corresponding construction in \cite{LL}.
    \label{no:cdr.1}
\end{nolabel}
\begin{nolabel}
    Let us fix notation first. Let $U$ be a differentiable manifold.
    Let $\cT$ be the tangent bundle of $U$ and $\cT^*$ be its
    cotangent bundle. We let $T=\Gamma(U, \cT)$ be the space of vector
    fields on $U$ and $A = \Gamma(U, \cT^*)$ be the space of
    differentiable $1$-forms on $U$. We let $\cC=\cC^\infty(U)$ be the
    space of differentiable functions on $U$. Denote by
    \begin{equation}
        <\,,\,> : A \otimes T \rightarrow \cC
        \label{cdr.2.1}
    \end{equation}
    the natural pairing. Finally, we denote by $\Pi$ the functor of
    \emph{change of parity}.

    Consider now a SUSY Lie conformal algebra $\cR$ generated by the
    super-vector space
    \begin{equation}
        \cC \oplus \Pi T \oplus A \oplus \Pi A
        \label{cdr.2.2}
    \end{equation}
    That is, we consider differentiable functions (to be denoted
    $f,g, \dots$) as even elements, vector fields $X, Y, \dots$ are
    odd elements, and finally we have two copies of the space of
    differential forms. For differential forms $\alpha, \beta, \dots
    \in A$
    we will denote the corresponding elements of $\Pi A$ by
    $\bar{\alpha}, \bar\beta, \dots$. The nonvanishing commutation
    relations in
    $\cR$ are given by (up to skew-symmetry):
    \begin{equation}
        \begin{aligned}
            {[X}_\Lambda f] &= X(f) \\
            [X_\Lambda Y] &= [X,Y]_{\mathrm{Lie}} \\
            [X_\Lambda \alpha] &= \lie_X \alpha + \lambda
            <\alpha, X> \\
            [X_\Lambda \bar{\alpha}] &= \overline{\lie_X \alpha} +
            \chi <\alpha, X>
        \end{aligned}
        \label{cdr.2.3}
    \end{equation}
    where $[,]_{\mathrm{Lie}}$ is the Lie bracket of vector fields
    and $\lie_X$ is the action of $X$ on the space of differential
    forms by the Lie derivative. The fact that (\ref{cdr.2.3}) is
    compatible with the Jacobi identity is a (long but) straightforward
    computation.

    We let $R(U)$ be the corresponding \emph{universal enveloping
    SUSY vertex algebra} of $\cR$. As noted in \cite{LL}, this
    algebra is too big. We want to impose some relations in $R(U)$.
    We let $1_U$ denote the constant function $1$ in $U$. Let $d: \cC
    \rightarrow A$ be the de Rham differential. Define $I(U) \subset
    R(U)$ to be the ideal generated by elements of the form:
    \begin{xalignat}{4}
        f_{(-1|1)}g - (fg), && f_{(-1|1)}X - (fX), &&
         f_{(-1|1)}\alpha - (f\alpha), &&  f_{(-1|1)} \bar{\alpha}
         - (\overline{f\alpha}), \\
         1_U - \vac, && Tf - df, && Sf - \bar{df} &&
        \label{eq:cdr.2.4}
    \end{xalignat}
    Finally we define the SUSY vertex algebra
    \begin{equation}
        \Omega^{\mathrm{ch}} (U) := R(U)/I(U)
        \label{eq:cdr.2.5}
    \end{equation}
    \label{no:cdr.2}
\end{nolabel}
\begin{thm}{\cite{LL}}\hfill
    \begin{enumerate}
    \item Let $M \subset \mathbb{R}^n$ be an open submanifold. The
    assignment $U \mapsto
    \Omega^{\mathrm{ch}}(U)$ defines a sheaf of SUSY vertex algebras
    $\Omega^{\mathrm{ch}}_M$ on $M$.
        \item For any diffeomorphism of open sets $M'
        \xrightarrow{\varphi} M $
         we obtain a canonical isomorphism of
        SUSY vertex algebras $\Omega^{\mathrm{ch}}(M)
        \xrightarrow{\Omega^{\mathrm{ch}}(\varphi)}
        \Omega^{\mathrm{ch}}(M')$. Moreover, given
        diffeomorphisms $M'' \xrightarrow{\varphi'} M'
        \xrightarrow{\varphi} M$, we have
        $\Omega^{\mathrm{ch}}(\varphi \circ \varphi') =
        \Omega^{\mathrm{ch}}(\varphi') \circ
        \Omega^{\mathrm{ch}}(\varphi)$.
    \end{enumerate}
    \label{thm:cdr.3}
\end{thm}
\begin{nolabel}
    This theorem allows us to construct a sheaf of SUSY vertex
    algebras
    in the Grothendieck topology on $\mathbb{R}^n$ (generated by open
    embeddings). This in turn lets us attach to any smooth manifold
    $M$, a sheaf of SUSY vertex algebras $\Omega^{\mathrm{ch}}_M$. We
    call this sheaf the \emph{chiral de Rham complex of $M$}.
    \label{no:cdr.3.b}
\end{nolabel}
\begin{rem}
    In the algebraic case, this construction gives the
    \emph{chiral de Rham complex} as described in the previous section in
    terms of coordinates. Indeed, we see that identifying
    $B^i$ with the field corresponding to the coordinate $x^i$ and
    $\Psi_i$ with the field cooresponding to the vector field
    $\partial_{x^i}$, the relations defining $I(U)$ are obviously
    satisfied.
    \label{rem:cdr.4}
\end{rem}

\section{Recollections on hyperk\"{a}hler manifolds} \label{hyperkaler}

In this section, we briefly review the notion of a hyperk\"{a}hler manifold following \cite{Joyce}. Let $(M,g)$ be a Riemannian manifold of real dimension $2n$, and $J$ a complex structure on $M$. The metric $g$ is Hermitian if $J$ is an isometry of $g$, i.e. if
\[
g(J u, J v) = g(u,v) \; \; \; \text{for tangent vectors} \; u, v
\]
Given $(M,g,J)$ we can define a $2$--form $\omega$ by:
\[
\omega (v, w) = g(J v, w)
\]
$\omega$ is called the \emph{Hermitian form} of $g$. The metric $g$ is said to be
\emph{K\"{a}hler} if $d \omega = 0$, in which case $\omega$ is called the \emph{K\"{a}hler
form}. The following theorem \cite{Joyce} gives other useful characterizations of K\"{a}hler metrics.

\begin{thm}
Let $(M,g,J)$ be as above, and $\nabla$ denote the Levi-Civita connection of $g$. Then the following conditions are equivalent:
\begin{enumerate}
\item $g$ is K\"{a}hler
\item $\nabla J = 0$
\item $\nabla \omega = 0$
\item The holonomy group of $g$ is contained in $U(n) \subset SO(2n, \mathbb{R})$
\end{enumerate}
\end{thm}

 A Riemannian manifold  $(M,g)$ of dimension $4n$ is \emph{hyperk\"{a}hler} if it possesses three complex structures $I, J, K$, such that $(M,g,I), (M,g,J),$ and $(M,g,K)$ are each K\"{a}hler, and satisfy the quaternionic relations
 \[
 IJ = -JI = K
 \]
 This is equivalent to the holonomy group of $g$ lying inside of $Sp(n,\mathbb{H}) \subset SO(4n,\mathbb{R})$. Dualizing using the metric yields three K\"{a}hler forms $\omega_I, \omega_J, \omega_K$.

\section{Superconformal structures from geometry} \label{geometry}

\begin{nolabel}
    Let now $(M,g)$ be a Riemannian manifold. Denote by $\Gamma^i_{j_k}$ the
    Christoffel symbols of the Levi-Civita connection of $M$. Let $I$ be an
    endomorphism of the tangent bundle of $M$, namely $I$ is a tensor of type
    $1,1$ on $M$. Let
    $\omega_i^j$ be the coordinate components of such a tensor. We will systematically raise and
    lower indices using the metric $g$ and sum over repeated indices.
    \label{no:3.4}
\end{nolabel}

In the rest of this paper, we will adopt the following notational convention, aimed at reducing clutter:
\bigskip

\noindent {\bf Convention:} In the $\beta \gamma - b c$ SUSY vertex algebra,
expressions such as
\[
 (\omega_{i}^{\; j} SB^i)\Psi_j +
        \Gamma^i_{jk}
        \omega_i^{\;j} TB^k
\]
correspond to the field
\begin{multline}
 (\omega_{i}^{\; j} (B^1(z,\theta), \cdots, B^n(z, \theta))
 (SB^i)(z,\theta))\Psi_j(z,\theta) + \\ +
        \Gamma^i_{jk} (B^1(z,\theta), \cdots, B^n(z, \theta))
        \omega_i^{\;j} (B^1(z,\theta), \cdots, B^n(z, \theta))
        (TB^k)(z,\theta)
\end{multline}
In other words, all coefficients of tensors, Christoffel symbols,
etc. are being evaluated on the superfields $B^i(z,\theta)$. This in
turn should be interpreted as follows. If $f(x^1, \cdots, x^n)$ is a
$C^{\infty}$ function on an open set $U$, then
\begin{align}
f(B^1(z,\theta), \cdots, B^n (z, \theta)) &= f(b^1(z) + \theta \phi^1 (z), \cdots, b^n (z) + \theta \phi^n (z)) \\
                                                            &= f(b^1(z), \cdots, b^n (z)) + \theta \sum^{n}_{i=1} \frac{\partial f}{\partial x^i} (b^1(z), \cdots, b^n (z)) \phi^i (z)
\end{align}
The meaning of expressions such as $f(b^1(z), \cdots, b^n (z))$ is
explained in \cite{malikov} (see also section \ref{Cinfty} and
\cite{LL}).
\begin{lem}
    The assignment:
    \begin{equation}
        I \mapsto J = (\omega_{i}^{\; j} SB^i)\Psi_j +
        \Gamma^i_{jk}
        \omega_i^{\;j} TB^k
        \label{eq:3.5.1}
    \end{equation}
    defines a linear morphism
    \begin{equation}
        \Gamma (M, \mathscr{E}nd(T_M)) \rightarrow \Gamma(M, \Omega^{\mathrm{ch}}_M)
        \label{eq:3.5.2}
    \end{equation}
    \label{lem:3.5}
\end{lem}
\begin{nolabel}
    In the case where $I$ is a complex structure on $M$ with associated K\"{a}hler
    form $\omega$, we will denote the corresponding current $J$ by $J_\omega$.
    \label{no:caca5.b}
\end{nolabel}
\begin{proof}
    Given that the fields $B^i$ transform as coordinates do, we will simplify
    the notation and denote:
    \begin{equation}
        \frac{\partial\tilde{B}^i}{\partial B^j} = \frac{\partial g^i
        (B)}{\partial B^j}.
        \label{eq:3.5.p.1}
    \end{equation}
    The first term in $J_\omega$ expressed in the coordinates $\tilde{B}^i$
    is given by:
    \begin{equation}
        \left( \frac{\partial B^k}{\partial \tilde{B}^i}
        \frac{ \partial \tilde{B}^j}{ \partial B^l} \frac{ \partial
        \tilde{B}^i}{\partial B^m} \omega_k^{\;l} SB^m \right) \left(
        \frac{\partial B^n}{\partial \tilde{B}^j} \Psi^n \right) = \left(
        \frac{ \partial \tilde{B}^j}{ \partial B^l}
         \omega_k^{\;l} SB^k \right) \left(
        \frac{\partial B^n}{\partial \tilde{B}^j} \Psi^n \right).
        \label{eq:3.5.p.2}
    \end{equation}
    Using quasi-associativity (\ref{eq:2.8.2}) we see that this is:
    \begin{equation}
        \left( \omega_{i}^{\; j} SB^i \right) \Psi_j -
        \frac{\partial^2 B^l}{\partial \tilde{B}^j \partial\tilde{B}^m}
        \frac{\partial \tilde{B}^j}{ \partial{B}^k} \frac{\partial
        \tilde{B}^m}{ \partial B^n} \omega_l^{\; k} TB^n
        \label{eq:3.5.p.3}
    \end{equation}
    On the other hand, in the second term in $J_\omega$, there are no
    quasi-associativity issues and the anomalous term comes from the
    transformation properties of the Christoffel symbols. Indeed, the second
    term in $J_\omega$ transforms as:
    \begin{equation}
        \frac{\partial^2 B^l}{\partial \tilde{B}^j \partial \tilde{B}^k}
        \frac{\partial \tilde{B}^i}{ \partial B^l} \frac{\partial
        B^m}{\partial \tilde{B}^i} \frac{\partial
        \tilde{B}^j}{\partial B_n} \frac{\partial \tilde{B}^k}{ \partial
        B^p} \omega_{m}^{\;n} TB^p + \Gamma^i_{jk} \omega_{i}^{\;j} TB^k
        \label{eq:3.5.p.4}
    \end{equation}
    Adding (\ref{eq:3.5.p.4}) and (\ref{eq:3.5.p.3}) we obtain the result.
\end{proof}
\begin{thm}\hfill
\begin{enumerate}
    \item   Let $(M,g)$ be a Riemannian manifold of dimension $n$, and
        $\textbf{g}=\log\sqrt{\det g_{ij}}$, where $\det g_{ij}$ is the
    determinant of the metric. Then
    \begin{equation}
        H = SB^i S\Psi_i + TB^i \Psi_i - TS \textbf{g}
        \label{eq:3.7.1}
    \end{equation}
    generates an $N=1$ superconformal structure of central charge $3 n$. We shall refer to $H$ as the Neveu-Schwarz vector.
    \item Let $(M,g)$ be a Calabi Yau $2 n$-manifold with
    K\"{a}hler form $\omega$. Then $J_\omega$ and $H$ generate an $N=2$ vertex algebra
    structure of central charge $6n$.
    \item If moreover $M$ is hyperk\"{a}hler of dimension $4n$, with three K\"{a}hler structures
        $\omega,\, \eta$ and $\gamma$ such that the corresponding complex
        structures satisfy the quaternionic relations, then $J_\omega, \,
        J_\eta, \, J_\gamma$ and $H$ generate an $N=4$ vertex algebra of
        central charge $12 n$.
\end{enumerate}
\label{thm:3.7}
\end{thm}
\begin{proof}\hfill
    \begin{enumerate}
        \item
    If we write $H = H^0 - TS \textbf{g}$, then the fact that $H^0$ defines a
    global section of $\Omega^{\mathrm{ch}}_M$ and that $H^0$ is a Neveu
    Schwarz vector of central charge $3 n$ follows from the analogous
    results in \cite{malikov}. The fact that $TS \textbf{g}$ is a well defined global
    section follows since $M$ is orientable. To check that $H$ is indeed a
    Neveu Schwarz vector (\ref{eq:2.9.2}) we compute:
    \begin{equation}
        {[}{H^0}_\Lambda \textbf{g}] = (2T + \chi S) \textbf{g}
        \label{eq:3.7.p.1a}
    \end{equation}
    therefore
    \begin{equation}
        \begin{aligned}
            {[}{H^0}_\Lambda TS \textbf{g}] &= (\lambda+ T) (S+\chi)
            (2T + \chi S) \textbf{g} \\
            &= (2T + 3\lambda + \chi S) TS \textbf{g} + \lambda^2 S
            \textbf{g} + \lambda \chi T \textbf{g} \\
            [TS \textbf{g}_\Lambda H^0] &= - \lambda \chi (2T - (\chi
            +S)S ) \textbf{g}\\
            &= - \lambda \chi T \textbf{g} - \lambda^2 S
            \textbf{g}
        \end{aligned}
        \label{eq:3.7.p.1b}
    \end{equation}
    It follows then that
    \begin{multline}
        [H_\Lambda H] = [{H^0}_\Lambda H^0] - [TS \textbf{g}_{\Lambda} H^0] -
        [{H^0}_\Lambda TS \textbf{g}] = \\ = (2T + 3\lambda + \chi S) H^0
        + n \lambda^2  \chi -
        (2T + 3\lambda + \chi S)TS \textbf{g} = \\ = (2T + 3\lambda + \chi S) H
        + n \lambda^2  \chi
        \label{eq:3.7.p.1c}
    \end{multline}

    \item
    To check the remaining commutation relations of the $N=2$ vertex algebra
    as in example \ref{ex:2.11.a}, it is enough to do it in any
    coordinate system. In particular we may choose holomorphic coordinates
    $\{x_\alpha\}$ (resp. anti-holomorphic coordinates
    $\{x_{\bar{\alpha}}\}$) for the complex structure associated to $\omega$ such that
    \begin{equation}
        \omega_i^{\; j} = \begin{pmatrix}
            i \Id & 0 \\ 0 & -i \Id
        \end{pmatrix}.
        \label{eq:3.7.p.1}
    \end{equation}
    In this case $J= J_\omega$ is given by
    \begin{equation}
        J= i SB^\alpha \Psi_\alpha - i SB^{\bar\alpha} \Psi_{\bar\alpha} + i
         \textbf{g}_{,\alpha} TB^\alpha - i \textbf{g}_{,\bar\alpha}
         TB^{\bar\alpha}.
        \label{eq:3.7.p.3}
    \end{equation}
    Here we have used the fact that for a K\"{a}hler manifold
    \begin{equation}
        \Gamma^\alpha_{\gamma \alpha} = \textbf{g}_{,\gamma} =
        \partial_\gamma \textbf{g} \qquad \Gamma^{\bar
        \alpha}_{\bar{\gamma} \bar{\alpha}} = \textbf{g}_{, \bar{\gamma}}
        = \partial_{\bar{\gamma}} \textbf{g}.
        \label{eq:3.7.p.4}
    \end{equation}
    Let us first compute $[H_\Lambda J]$. For this we need
    \begin{equation}
        \begin{aligned}
            {[}H_\Lambda \Psi_\alpha] &= (2T + \lambda + \chi S)
            \Psi_\alpha +
        \lambda \chi \textbf{g}_{, \alpha} \\ [H_\Lambda
        \Psi_{\bar{\alpha}}] &= (2T + \lambda + \chi S)
        \Psi_{\bar{\alpha}} + \lambda \chi \textbf{g}_{,\bar{\alpha}} \\
        {[}H_\Lambda B^\alpha ] &= (2T + \chi S) B^\alpha \\
        [H_\Lambda B^{\bar{\alpha}}] &= (2 T + \chi S) B^{\bar{\alpha}} \\
        [H_\Lambda SB^\alpha] &= (2T + \lambda + \chi S) SB^\alpha \\
        [H_\Lambda SB^{\bar{\alpha}}] &= (2T + \lambda + \chi S)
        SB^{\bar{\alpha}}
    \end{aligned}
        \label{eq:3.7.p.5}
    \end{equation}
    Using now the non-commutative Wick formula we obtain
    \begin{multline}
        [H_\Lambda SB^\alpha \Psi_\alpha] = \left( (2T + \lambda + \chi
        S) SB^\alpha \right) \Psi_\alpha + SB^\alpha (2T + \lambda + \chi
        S) \Psi_\alpha + \\ + SB^\alpha \lambda \chi \textbf{g}_{,\alpha}
        + \int_0^\Lambda [ (2T + \lambda + \chi S) {SB^\alpha}_\Gamma
        \Psi_\alpha] d \Gamma = \\
        (2 T + 2 \lambda + \chi S) SB^\alpha \Psi_\alpha -\lambda \chi
        \textbf{g}_{,\alpha} SB^\alpha + \int_0^\Lambda (-2 \gamma + \lambda - \chi \eta)
        \eta [B^\alpha_\Lambda \Psi_\alpha] d\Gamma
        \label{eq:3.7.p.6}
    \end{multline}
    Since the integral
    clearly vanishes, we obtain:
    \begin{equation}
        [H_\Lambda SB^\alpha \Psi_\alpha] = (2T + 2 \lambda +\chi S)
        SB^\alpha \Psi_\alpha - \lambda \chi \textbf{g}_{,\alpha}
        SB^\alpha.
        \label{eq:3.7.p.7}
    \end{equation}
    Similarly, we compute now
    \begin{equation}
        \begin{aligned}
            {[\textbf{g}_{,\alpha}}_\Lambda H] &= SB^i (\chi + S)
            \textbf{g}_{, \alpha i} + TB^i \textbf{g}_{,\alpha
            i} \\
            &= - \chi S \textbf{g}_{,\alpha} + SB^i
            S\textbf{g}_{,\alpha i} + T \textbf{g}_{,\alpha} \\
            &= (T - \chi S) \textbf{g}_{,\alpha} -
            \textbf{g}_{,\alpha i j} SB^j SB^i \\
            &= (T- \chi S) \textbf{g}_{,\alpha} \\
            [H_\Lambda {\textbf{g}_{,\alpha}}] &= (2T + \chi S)
            \textbf{g}_{,\alpha}
        \end{aligned}
        \label{eq:3.7.p.8}
    \end{equation}
     We also have:
    \begin{equation}
        [H_\Lambda TB^\alpha] = (\lambda + T) (2T + \chi S) B^\alpha =
        (2T + 2\lambda + \chi S)T B^\alpha + \lambda \chi S B^\alpha
        \label{eq:3.7.p.9}
    \end{equation}
    Hence using the Wick formula again and noting that the integral term
    trivially vanishes, we obtain:
    \begin{multline}
        [H_\Lambda \textbf{g}_{,\alpha} TB^\alpha] = \left((2T + \chi S)
        \textbf{g}_{,\alpha} \right) TB^\alpha + \\
        + \textbf{g}_{,\alpha} (2T
        + 2\lambda + \chi S) TB^\alpha + \textbf{g}_{,\alpha} \lambda
        \chi SB^\alpha = \\ = (2T + 2\lambda +\chi S)
        \textbf{g}_{,\alpha} TB^\alpha + \lambda \chi
        \textbf{g}_{,\alpha} SB^\alpha
        \label{eq:3.7.p.10}
    \end{multline}
    Adding (\ref{eq:3.7.p.7}) to (\ref{eq:3.7.p.10}) plus their conjugates we
    obtain:
    \begin{equation}
        {[}H_\Lambda J] = (2 T + 2\lambda + \chi S) J
        \label{eq:3.7.p.11}
    \end{equation}
    as we wanted.

    Finally, we need to check (\ref{eq:2.11.a.2}). For this we compute:
    \begin{equation}
        \begin{aligned}
            {[\Psi_\alpha}_\Lambda J] &= i \chi \Psi_\alpha + i
            \textbf{g}_{,\beta \alpha} TB^\beta + i \lambda
            \textbf{g}_{,\alpha} \\
            &= i \chi \Psi_\alpha +  i (T + \lambda)
            \textbf{g}_{,\alpha} \\
            [J_\Lambda \Psi_{\alpha}] &= - i (\chi + S) \Psi_\alpha - i
            \lambda \textbf{g}_{,\alpha} \\
            {[SB^\alpha}_\Lambda J] &= - i\chi SB^\alpha \\
            [J_\Lambda SB^\alpha ] &= i(\chi + S)SB^\alpha \\
        \end{aligned}
        \label{eq:3.7.p.12}
    \end{equation}
    Here in the second line we used the fact that $M$ is Ricci flat,
    therefore $\textbf{g}_{,\alpha \bar{\beta}} = 0$. Using this we can compute now:
    \begin{multline}
        [J_\Lambda SB^\alpha \Psi_\alpha] = i \left( (\chi+ S) SB^\alpha
        \right) \Psi_\alpha +\\+ i SB^\alpha (\chi + S) \Psi_\alpha + i
        SB^\alpha \lambda \textbf{g}_{,\alpha} + \int_0^\Lambda i
        [{(\chi+S)SB^\alpha}_\Gamma \Psi_\alpha] d\Gamma = \\
        i TB^\alpha \Psi_\alpha + i SB^\alpha S\Psi_\alpha + i\lambda
        SB^\alpha \textbf{g}_{,\alpha} + i n \int_0^\Lambda (\eta -
        \chi)\eta  d\Gamma   \\ = i TB^\alpha \Psi_\alpha + i SB^\alpha
        S\Psi_\alpha  + i \lambda SB^\alpha \textbf{g}_{,\alpha} + i n
        \lambda \chi
        \label{eq:3.7.p.13}
    \end{multline}
    Similarly we have:
    \begin{equation}
        \begin{aligned}
            {[\textbf{g}_{,\alpha}}_\Lambda J] &= - i SB^\beta
            \textbf{g}_{,\alpha \beta} \\
            &= - i S \textbf{g}_{,\alpha} \\
            [J_\Lambda \textbf{g}_{,\alpha}] &= - i S
            \textbf{g}_{,\alpha} \\
            [J_\Lambda TB^\alpha] &= - i (\lambda + T) SB^\alpha
        \end{aligned}
        \label{eq:3.7.p.14}
    \end{equation}
    Hence using the Wick formula we obtain:
    \begin{multline}
        [J_\Lambda \textbf{g}_{,\alpha} TB^\alpha] = - i (S
        \textbf{g}_{,\alpha} )TB^\alpha - i \textbf{g}_{,\alpha} (\lambda
        + T) SB^\alpha = \\ =  - i \lambda \textbf{g}_{,\alpha} SB^\alpha
        - i S(\textbf{g}_{,\alpha} TB^\alpha)
        \label{eq:3.7.p.15}
    \end{multline}
    Adding (\ref{eq:3.7.p.15}) and (\ref{eq:3.7.p.13}) plus their conjugates,
    we obtain:
    \begin{equation}
        \begin{aligned}
            {[J}_\Lambda J] &= - TB^i \Psi_i - SB^i S\Psi_i - \lambda \chi 2 n +
        S (\textbf{g}_{,\alpha} TB^\alpha) + S
        (\textbf{g}_{,\bar{\alpha}} TB^{\bar{\alpha}}) \\
        &= - TB^i \Psi_i - SB^i S\Psi_i + ST \textbf{g} - 2n \lambda \chi
    \end{aligned}
        \label{eq:3.7.p.16}
    \end{equation}
    \item In order to prove the rest of the statement, we may choose
        holomorphic coordinates for one of the three complex structures.
        We suppose then that $J_\omega$ is given by (\ref{eq:3.7.p.3}).
        This forces the other two complex structures to be of the form:
        \begin{equation}\eta_{i}^{\; j} =
            \begin{pmatrix}
                0 & \eta_\alpha^{\; \bar\beta} \\
                \eta_{\bar \alpha}^{\; \beta} & 0
            \end{pmatrix}
            \qquad \gamma_i^{\; j} =\begin{pmatrix}
                0 & i \eta_\alpha^{\; \bar \beta} \\
                - i \eta_{\bar \alpha}^{\; \beta} & 0
            \end{pmatrix}.
            \label{eq:3.7.p.17}
        \end{equation}
        with $\eta_\alpha^{\;\bar \beta} \eta_{\bar \beta}^{\; \gamma} =
        - \delta_{\alpha}^{\; \gamma}$.  It follows then that
        $J^+ = \tfrac{1}{2} (J_\eta - i J_\gamma)$ is given by
        \begin{equation}
            J^+ = \eta_\alpha^{\; \bar \beta} SB^\alpha \Psi_{\bar
            \beta}
            \label{eq:3.7.p.18}
        \end{equation}
        where we note that the second term in (\ref{eq:3.5.1}) vanishes
        in this case since $M$ is K\"{a}hler, therefore the only
        non-vanishing components of the Christoffel symbols are
        $\Gamma^\alpha_{\beta \gamma}$ and
        $\Gamma^{\bar{\alpha}}_{\bar{\beta}\bar{\gamma}}$. The theorem
        will be proved then if we show that\footnote{Note that the commutation relations with $J^- =
        \tfrac{1}{2}(J_\eta+iJ_\gamma)$ follow by
        conjugation.}
        \begin{equation}
            {[J_\omega}_\Lambda J^+] = i \left( 2 \chi + S \right) J^+.
            \label{eq:3.7.p.19}
        \end{equation}
        In order to do so we compute (we denote $J= J_\omega$):
        \begin{equation}
            \begin{aligned}
                {[\eta_\alpha^{\; \bar\beta}}_\Lambda J] &= - i
                {\eta_{\alpha}^{\; \bar\beta}}_{,\gamma} SB^\gamma +
                i
                {\eta_{\alpha}^{\;\bar{\beta}}}_{,\bar{\gamma}}
                SB^{\bar \gamma}  \\
                [J_\Lambda \eta_\alpha^{\;\bar\beta}] &=  -i
                {\eta_{\alpha}^{\; \bar\beta}}_{,\gamma} SB^\gamma +i
                {\eta_{\alpha}^{\;\bar{\beta}}}_{,\bar{\gamma}}
                SB^{\bar \gamma}
            \end{aligned}
            \label{eq:3.7.p.20}
        \end{equation}
        It follows then from (\ref{eq:3.7.p.12}) and the Wick formula
        \begin{equation}
            [J_\Lambda \eta_{\alpha}^{\; \bar\beta} SB^\alpha] = - i{
            \eta_\alpha^{\; \bar\beta}}_{,\gamma} SB^\gamma SB^\alpha
             + i {\eta_\alpha^{\; \bar\beta}}_{,\bar\gamma}
             SB^{\bar\gamma} SB^\alpha  + i \eta_\alpha^{\;
             \bar\beta} (\chi + S) SB^\alpha
            \label{eq:3.7.p.20a}
        \end{equation}
        Since the complex structure is parallel we have:
        \begin{equation}
            {\eta_{\alpha}^{\;\bar{\beta}}}_{,\gamma} = \Gamma_{\gamma
            \alpha}^\delta \eta_\delta^{\; \bar\beta} =
            {\eta_\gamma^{\;\bar{\beta}}}_{,\alpha}
            \label{eq:3.7.p.20b}
        \end{equation}
        Therefore the first term in (\ref{eq:3.7.p.20a}) vanishes and we
        have:
        \begin{equation}
            [J_\Lambda \eta_{\alpha}^{\; \bar\beta} SB^\alpha] =  i {\eta_\alpha^{\; \bar\beta}}_{,\bar\gamma}
             SB^{\bar\gamma} SB^\alpha  + i \eta_\alpha^{\;
             \bar\beta} (\chi + S) SB^\alpha
            \label{eq:3.7.p.20c}
        \end{equation}
        Now conjugating (\ref{eq:3.7.p.12}) we see that
        \begin{equation}
            [J_\Lambda \Psi_{\bar \beta}] = i (\chi + S)
            \Psi_{\bar\beta} + i \lambda \textbf{g}_{, \bar{\beta}}
            \label{eq:3.7.p.21}
        \end{equation}
        Now using the non-commutative Wick formula we obtain:
        \begin{multline}
            {[J_\omega}_\Lambda J^+] =
                     i \left(
                {\eta_{\alpha}^{\;\bar{\beta}}}_{,\bar{\gamma}}
                SB^{\bar \gamma}SB^\alpha \right) \Psi_{\bar \beta} + i
                \left( \eta_\alpha^{\; \bar \beta} (\chi + S)
                SB^\alpha \right) \Psi_{\bar\beta} - \\ - i \left(
                \eta_\alpha^{\;\bar\beta} SB^\alpha
                \right) (\chi + S) \Psi_{\bar \beta} - i \left(
                \eta_\alpha^{\;\bar\beta} SB^\alpha \right)\lambda
                \textbf{g}_{,\bar\beta} + \\ + i \int_0^\Lambda {[{\eta_{\alpha}^{\;\bar{\beta}}}_{,\bar{\gamma}}
                SB^{\bar \gamma}SB^\alpha}_\Lambda \Psi_{\bar\beta}] d
                \Gamma + i \int_0^\Lambda {[\eta_\alpha^{\; \bar \beta} (\chi + S)
                SB^\alpha}_\Lambda  \Psi_{\bar\beta}] d\Gamma
            \label{eq:3.7.p.22}
        \end{multline}
        Let us compute the integral term first. Clearly the second
        integral vanishes since $\Psi_{\bar\beta}$ commutes with
        $SB^\alpha$. The first term  on the other hand is given by:
        \begin{equation}
            - \int_0^\Lambda  \eta \;
            {\eta_\alpha^{\;\bar\beta}}_{,\bar\beta} SB^\alpha
            d\Gamma = - \lambda
            {\eta_\alpha^{\;\bar\beta}}_{,\bar\beta} SB^\alpha
            \label{eq:3.7.p.23}
        \end{equation}
        Re-grouping the terms in (\ref{eq:3.7.p.22}) we obtain:
        \begin{equation}
            {[J_\omega}_\Lambda J^+ ] = i (2 \chi + S) J^+ - i \lambda
            \left( \eta_{\alpha}^{\; \bar \beta}
            \textbf{g}_{,\bar\beta} SB^\alpha + {\eta_{\alpha}^{\;
            \bar \beta}}_{,\bar\beta}SB^\alpha\right)
            \label{eq:3.7.p.24}
        \end{equation}
        Finally, using again the fact that $\eta$ is parallel, we see
        immediately that:
        \begin{equation}
            {\eta_\alpha^{\;\bar \beta}}_{,\bar\beta} = -
            \Gamma^{\bar\beta}_{\bar \beta \bar \gamma}
            \eta_\alpha^{\; \bar \gamma} = - \textbf{g}_{,\bar\gamma}
            \eta_\alpha^{\; \bar \gamma}
            \label{eq:3.7.p.25}
        \end{equation}
        This proves (\ref{eq:3.7.p.19}) and the theorem.
    \end{enumerate}
\end{proof}
\begin{rem} \hfill
    \begin{enumerate}
        \item  In the Calabi-Yau case, choosing
            holomorphic coordinates $x^\alpha$ and
            antiholomorphic coordinates $x^{\bar{\alpha}}$,
            we see that $\textbf{g} = \textbf{g}_0 +
            \bar{\textbf{g}_0}$ where $\textbf{g}_0$ is
            holomorphic.
            The super-field $H$ can be decomposed as a sum of
            two terms $H_0 + \bar{H_0}$, where
            \begin{equation}
                H_0 = SB^\alpha S\Psi_\alpha + TB^\alpha
                \Psi_\alpha - TS \textbf{g}_0
                \label{rem.add.1}
            \end{equation}
            The super-field $J$ decomposes in a similar
            way as a ``holomorphic'' part $J_0$ and an
            ``anti-holomorphic'' part $\bar{J_0}$.
            These fields are invariant under holomorphic
            changes of coordinates, hence we obtain two
            commuting $N=2$ super conformal structures.

            We note that these fields are different from the
            ones considered in \cite{malikov}. In the case of
            the Virasoro field, the correction given by $TS
            \mathbf{g_0}$ appeared in \cite[page 16]{Witten}.
            When the metric is flat (i.e. $\mathbf{g_0} = 0$)
            we obtain the same topological structure as in
            \cite{malikov}.
        \item   When the manifold $M$ is complex but not
            Calabi-Yau, the decomposition $H = H_0 +
            \bar{H_0}$ is not invariant under holomorphic
            changes of coordinates. Therefore our $N=1$
            structure pairs in a non-trivial way the
            ``holomorphic'' and ``anti-holomorphic'' parts of
            the chiral de Rham complex of $M$
        \item   The fields $H$ and $J$ are defined for any almost complex manifold $M$ (though they
        generate $N=2$ only when $M$ is Calabi--Yau). In
            particular, the field $J$ allows us to construct
            a \emph{Dolbeault resolution} of the
            \emph{holomorphic chiral de Rham complex} in
            terms of the differentiable one (see also
            \cite{kapustin}). We plan to return to this
            matter elsewhere.
    \end{enumerate}
       \label{rem:3.8}
\end{rem}

\bibliographystyle{alpha}
\bibliography{refs2}

\end{document}